
\magnification=1200
\hsize=15.2truecm 
\vsize=21truecm
\voffset=1.2truecm 
\newdimen\hov
\newdimen\htol
\hov=0.3truecm 
\htol=0.05truecm 
\parskip \hov plus \htol minus \htol
\everypar{\noindent}

\font\bb=cmbx10

\font\eightbrom=cmbx8
\font\eightbita=cmmib10 scaled 500
\font\eightbsym=cmbsy10 scaled 500
\font\fifteenbrom=cmbx10 scaled 1439
\font\fifteenbita=cmmib10 scaled 1439
\font\fifteenbsym=cmbsy10 scaled 1439
\font\twbrom=cmbx10 scaled 1200
\font\twbita=cmmib10 scaled 1200
\font\twbsym=cmbsy10 scaled 1200

\def\mathpar {
 \textfont0=\fifteenbrom \scriptfont0=\twbrom
 \scriptscriptfont0=\eightbrom
 \textfont1=\fifteenbita \scriptfont1=\twbita
 \scriptscriptfont1=\eightbita
 \textfont2=\fifteenbsym \scriptfont2=\twbsym
 \scriptscriptfont2=\eightbsym
 \textfont3=\tenex \scriptfont3=\tenex}


\def\paragrafo #1
{{\fifteenbrom #1} \par \hfill\break\noindent }
\def\teo #1 {\vskip0.4truecm \noindent {\bb Theorem~#1~~~}}
\def\lem #1 {\vskip0.4truecm \noindent {\bb Lemma~#1~~~}}
\def\prop #1 {\vskip0.4truecm \noindent {\bb Proposition~#1~~~}}
\def\dff #1 {\vskip0.4truecm \noindent {\bb Definition~#1~~~}}
\def\oss #1 {\vskip0.4truecm \noindent {\bb Remark~#1~~~}}
\def\coro #1 {\vskip0.4truecm \noindent {\bb Corollary~#1~~~}}
\def\pdim  {\parskip 0pt \hfill\break \par
\noindent {\it Proof} \rm \hfill\break }
\def\cite #1 {\item{\bb [#1]~}}

\def\acapo {\hfill\break}

\def\hr {$\hookrightarrow$}
\def\lo{\longrightarrow} 

\def\pp2 {${\rm I}\!{\rm P}^{ 2}$}
\def\pn{${\rm I}\!{\rm P}^{ n}$}
\def\pN {${\rm I}\!{\rm P}^{ N}$}
\def\p5 {${\rm I}\!{\rm P}^{ 5}$}

\def\od {${\cal O} _{{\rm I}\!{\rm P}^{ n}} (d)$}
\def\on {${\cal O} _{{\rm I}\!{\rm P}^{ n}}$}
\def\fd  {$\varphi _{{\cal O} _{{\rm I}\!{\rm P}^{ n}}(d)}$}
\def\fm  {$\varphi _{{\cal O} _{{\rm I}\!{\rm P}^{ m}}(d)}$}
\def\f2  {$\varphi _{{\cal O} _{{\rm I}\!{\rm P}^{ 2}}(d)}$}
\def\fpiano2  {$\varphi _{{\cal O} _{{\rm I}\!{\rm P}^{ 2}}(2)}$}

\def\ho #1 {$H^0 ({\cal O} _ {{\rm I}\!{\rm P}^{ n}} (#1))$}

\def\e {$ E_d$}


{\fifteenbrom
\centerline {Syzygies of Veronese embeddings}
}
\vskip1truecm
{\centerline {Giorgio Ottaviani ~~and ~~Raffaella Paoletti\footnote{*}
{both authors were partially supported by GNSAGA and by Progetto Strategico 
"Calcolo Simbolico" of CNR.}}}
\vskip2truecm

\noindent
{\bf Abstract }~- ~{\it {We prove that the Veronese embedding  \fd $\, : \, $ \pn  $\, \hookrightarrow \, $ \pN ~ with $n\ge 2$, $d\ge 3$
does not satisfy property $N_p$ (according to Green and Lazarsfeld)
if $p\ge 3d-2$. We make the conjecture that also the converse holds. 
This is true for $n=2$ and for $n=d=3$.}}

\vskip2truecm
\paragrafo {Introduction}

Let \pn ~be the projective $n$-space over an algebraically closed 
field of characteristic zero 
and let \fd $\, : \,$\pn $\, $ \hr $\, $\pN ~be the 
Veronese embedding associated to the complete linear system $| $
\od $ | \, $. ~In order to understand the homogeneous ideal ${\cal I}
$~ of \pn ~in \pN ~ as well as its syzygies, it is useful to study some
properties about a minimal free resolution of ${\cal I} \,.$\acapo
M. Green and R. Lazarsfeld (see [G2], [GL]) introduced the property $N_p$
(Definition $1.3 )$~ for a complete projective non singular variety ~ X 
$\, $ \hr $\, $\pN  ~~embedded in \pN ~ with an ample line bundle 
$L \,.$~ When property $N_p$~ holds for every integer $p$ the
resolution of $\cal I$ ~is "as nice as possible".~ M. Green proves in [G2], 
Theorem 2.2, that \fd ~ satisfies $N_p$ if $p\le d\, $. L. Manivel has generalized 
this result to flag manifolds [M]. The rational
normal curves (which are the Veronese embeddings of ${\rm I}\!{\rm P} ^1
)$~ satisfy $N_p ~~\forall p \,.$~ C. Ciliberto showed us that the results 
of [G1] imply that \f2 ~with $d\geq 3 $ satisfies $N_p$
 if
$p\leq 3d-3 \,.$ ~ This sufficient condition has been found also by 
C. Birkenhake  in [B1] as a corollary of a more general result.
Here we prove that this  condition is also necessary (Theorem $3.1)$~ and we formulate 
 (for $n\ge 2$) the following
\vskip1truecm
\noindent
{\bf Conjecture  }~~~ 
 \fd ~~  {\it satisfies }$N_p ~~ \Longleftrightarrow ~~ \left\{\matrix{
 n=2, d=2, ~\forall p \cr n\ge 3, d=2, ~p\le 5\cr n\ge 2, d\ge 3, ~p\le 3d-3}\right . ~. $
\medskip
Our precise result is the following:

\teo
{\it {The implication '$\Longrightarrow$' of the previous conjecture is true.}}
\medskip
Moreover we remark that the implication
'$\Longleftarrow$' of the previous conjecture is true
in the cases $n=2$ ([G1]), $n=d=3$ ([G1]), $d=2$ ([JPW]). This solves the Problem 
4.5 of [EL] (raised by Fulton) in the first cases given by the projective plane
and by the cubic embedding of the projective 3-dimensional space.

We remark also that our conjecture would be overcomed by the knowledge of the 
minimal resolution of the Veronese variety. This is stated as an open problem in 
[G2] (remark of section 2).  Our results can be seen as a step towards this problem.

        The paper is organized as follows: in {\it
section} $1$~ we recall some definitions we will need later and we
improve a known cohomological criterion for the property $N_p$; in {\it
section} $2$~ we prove our main results
and in {\it section} $3$ ~we fit our results into the literature.

We thank the referee for the improvements suggested in the presentation of the paper.

\acapo
\indent SUBJECT CLASSIFICATION: 14M20, 14F17. 

\vskip2truecm

\paragrafo{ 1 - Notations and preliminaries}

Let $V$~ be a vector space of dimension $n+1$~over an
algebraically closed field $I \! \! K $~of characteristic 0 and let \pn $\, =
\, {\rm I}\!{\rm P} (V ^* )$~the projective space associated
to the dual space of $V \,.$ \acapo
Note that $ 
H^0 \left( {\rm I}\!{\rm P}^n , \; {\cal O}_{{\rm I}\!{\rm P}^n}
(d) \, \right)  \, \cong \, S^d V ~~~ \forall d \geq 0\,.$ 
\acapo
For any vector bundle $E~{\rm over} ~{\rm I}\!{\rm P}^n \, ,$~
we will denote by $H^i (E \, )$~ the $i$-th cohomology group
of $E ~ {\rm over} ~{\rm I}\!{\rm P}^n $~and by $E(t)$~the
tensor product $E \, \otimes \, {\cal O}_{{\rm I}\!{\rm P}^n}
(t) \, .$

The following bundles will play a fundamental role in this
paper:

\dff{1.1} For any positive integer $d \, ,$~the line bundle
\on $(d)$~is gene-\break
rated by global sections ~\ho {d} 
$\, \cong \, S^d V$~so that the evaluation map\break 
$ ev \, : \; S^d V
\otimes $\on $\, \lo \, $\on $(d)$~is surjective. Call \e ~the kernel.
Thus the vector bundle \e ~is defined by the exact sequence:
$$0 ~\lo ~ E_d
~\lo  ~S^{d\,} V\, \otimes \, 
{\cal O}_{ {\rm I}\!{\rm P}^n}
~~~ {
\setbox2=\hbox{$-\hskip -5pt -\hskip -5pt 
\longrightarrow$}
\setbox1=\hbox to \wd2{\hfill $\scriptstyle ev$  \hfill}
\hbox {\hskip -8pt \vbox {\box1 \vskip -7pt \box2}}
\hskip 2pt}    ~~
{\cal O}_{ {\rm I}\!{\rm P}^n} \,(d)~~
\lo ~  0 \, . \leqno (1.2)$$

It follows immediately from the definition that the bundle \e
~has rank\break
$N := rk_{E_d} \, = \, {n+d \choose n} -1 $~and first Chern class
$c_1 (E_d ) = -d \,.$ 
\vskip0.1truecm 

Note that, if $d\,=\,1\, ,~ (1.2)$~is the
dualized Euler sequence so that: $$ E_1 ~\cong ~ \Omega ^1_{
{\rm I}\!{\rm P}^n}(1) ~~~{\rm and }~~
\bigwedge ^q E_1 \, \cong \, \, \Omega ^q _{{\rm I}\!{\rm P}^n}(q) \, .$$ 

For any integer $d \geq 0 \, ,$~ we will denote by \fd ~the {\it 
Veronese ~embedding}~ \fd $\, : \, $ \pn $\, $ \hr $\, $ \pN ~
associated to the complete linear system
$|$ \on $(d) \, | $~ of dimension $N+1 \, := \,  {n+d \choose n} \,.$~ 
Recall that if  
 $[x_0 : \, \dots \, : x_n ]$~ is a system of homogeneous
coordinates on \pn ~and 
$[y_0  : \, \dots \, : y_N ]$~ on \pN $\, = \, 
{\rm I}\!{\rm P} \left( H^0 \left( {\cal O}_{{\rm I}\!{\rm P}^n}
(d) \right) ^* \right) \, ,$~then   \fd ~is the embedding:
$$
[x_0 : \, \dots \, : x_n ] \,  \hookrightarrow \, 
[x_0^d : x_0^{d-1}x_1 : \, \dots \, : x_n^d ] \,.$$

\bigskip
With the above notation, let $S \, := \, \bigoplus _{k \geq 0} \,
S^k \, \left( H^0 \left( {\rm I}\!{\rm P}^n , \; {\cal O}_{{\rm I}\!{\rm
P}^n}(d) \right) \right) $~ be the homogeneous 
coordinate ring of \pN
~ and define the graded $S$-module\break
$R \, := \,  \bigoplus _{k \geq 0}
\,  H^0 \left( {\rm I}\!{\rm P}^n , \; {\cal O}_{{\rm I}\!{\rm P}^n}(kd)
\right) \, .$~ Let 
$$ 0 \,
\lo \,  \oplus _j \, S (-j) ^{{\widetilde b}_{rj}} 
\, \lo \, \cdots \, \lo \,
\oplus _j \, S (-j) ^{{\widetilde b}_{0j}} \, \lo \, R \, 
\lo \, 0 $$
be a minimal free resolution of $R$~with {\it graded Betti numbers}~$
{{\widetilde b}_{ij}} \,.$ 

\bigskip

\dff {1.3} For any integer $p \geq 0 $~ the embedding 
\fd $\, : \, $ \pn $\, $ \hr $\, $ \pN ~ is said to satisfy {\it
property}~ $N_p$~ if 
$$ 
{{\widetilde b}_{0j}} = \cases { 1 \cr 0 \cr } ~{\rm if} ~\eqalign 
{j=0 {\phantom {wise}}\cr {\rm otherwise } \cr } ~~~~~{\rm and}~~~~~~~ 
{{\widetilde b}_{ij}} =
0 ~~~~~{\rm for} ~~j\not= i+1 ~~{\rm when }~~1 \leq i \leq p \,.$$
Thus: \acapo
\indent
$N_0$ ~ means that \fd (\pn ) ~is projectively normal in \pN $\, ;$\acapo
\indent 
$N_1$ ~ means that $N_0$ ~holds and the ideal $I$~ of \pn $\, $\hr $\, $\pN
~is generated  by quadrics;\acapo
\indent 
$N_2$ ~ means that, moreover, the module of syzygies among quadratic
generators $Q_i \in I $~ is spanned by the relations of the form $\sum L_i
Q_i \, = \, 0$~ where the $L_i$~are {\it linear } polynomials; \acapo
\indent
and so on.

 \bigskip

\oss {1.4} Let ${\cal C} \, \hookrightarrow \, {\rm I}\!{\rm P}^d$~ be 
the rational normal curve (of degree $d$) in ${\rm I}\!{\rm P}^d
\,.$\break 
If $V$~is a vector space of dimension $2\,,$~then 
 ${\cal C} \cong {\rm I}\!{\rm P}(V^*) \, \hookrightarrow \, 
{\rm I}\!{\rm P}^d = {\rm I}\!{\rm P}(S^d V ^{\, *}) $~~is the image of the
Veronese embedding $\varphi _{{\cal O}_{{\rm I}\!{\rm P}^1 (d)}} \, : \, 
{\rm I}\!{\rm P}^1 \, \hookrightarrow \, {\rm I}\!{\rm P}^d \,.$ 

It is well known (e.g. by using the Eagon-Northcott complex)
that the sheaf ideal ${\cal I}$~of ${\cal C}$~ in 
${\cal O}_{{\rm I}\!{\rm P}^d}$~ has the following resolution:
$$
0  \lo  
{\cal O}_{{\rm I}\!{\rm P}^d} (-d) ^{\oplus b_d}  \lo
{\cal O}_{{\rm I}\!{\rm P}^d} (-d+1) ^{\oplus b_{d-1}}  \lo 
   \cdots  \lo 
 {\cal O}_{{\rm I}\!{\rm P}^d} (-2) ^{\oplus b_2} 
 \lo   {\cal I}   \lo   0 $$
where $b_k \, := \, (k-1) {d \choose k} \,$. So the Veronese embeddings of 
${\rm I}\!{\rm P}^1 $~ satisfy $N_p ~~\forall p \,.$\acapo

From [B2] , Remark 2.7, and [G1] we have the following cohomological criterion:

\prop{1.5} {\it {The Veronese embedding ~ \fd ~~satisfies property
$N_p$~if and only if }}
$$ H^1 \left( \bigwedge ^q  E_d \, (jd)
\, \right) \, = \, 0 \hskip 2truecm {\it { for ~~1\leq q \leq
p+1 ~~ and}} ~~\forall j\geq 1 \,.$$
$\hfill \diamond $
\vskip1truecm

We have the following cohomological criterion, which improves slightly the previous one (in fact $H^2  \left( \bigwedge^{q}  E_d \right) \simeq H^1 \left( \bigwedge^{q-1} E_d(d) \right) \, )$.

\teo {1.6}  
{\it {The Veronese embedding ~\fd ~~  satisfies property ~ $N_p $~  if and only if}}  
~ $H^2 \left( \bigwedge^{q}  E_d \right) =0  ~~{\it for}~ 1\leq q\le p+2 \,.$

\bigskip
The proof of Theorem 1.6 relies on the following

 \prop {1.7} 
 \hfill {\it If } ~~$H^2( \bigwedge^q E_d)=0 ~~{\it for }~~ 1\le q\le k \, ,~~ 
{\it then } ~~H^2( \, \bigwedge^q E_d(t))=0\break 
{\it for }~~ 1\le q\le k ~~{\it and} ~~\forall t\ge 0 \,.$
\pdim 
Consider the two exact sequences:
\medskip
\noindent 
$ (*)  \hskip2truecm 0 \lo \bigwedge^q E_d (t-1) \lo \bigwedge^q E_d (t) \lo  \bigwedge^q E_d (t)|_{{\rm I}\!{\rm P}^{n-1}}
 \lo 0 \, , \hfill $
\medskip
\noindent 
$ (**) \hskip 0.2truecm 0 \lo \bigwedge^q E_d (t-1) \lo
 \bigwedge^q (S^dV)\otimes
{\cal O}_{{\rm I}\!{\rm P}^n}(t-1) \lo
\bigwedge^{q-1} E_d (t+d-1)\lo 0  \, . \hfill $
\medskip
\noindent
The proof is by double induction on $n$ and $k$. The statement is true for 
$n=2$ ~(Serre duality) 
and for $k=1$ ~(it follows immediately from $(1.2) $). 
 From the cohomology sequence associated to ~$(**)$ ~ with $t=0$ ~and 
the inductive hypothesis on $k$ ~we get  $H^3 (\bigwedge^q E_d (-1) )=0 
~{\rm for }~ 1\le q\le k\, .$ ~ Since $E_d |_{{\rm I}\!{\rm P}^{n-1}} \,
\cong \, {\widetilde 
E_d} \oplus {\cal O} ^{\oplus {n+d-1 \choose n}} _{{\rm I}\!{\rm P}^{n-1}}  , ~ {\rm where} ~{\widetilde E_d} $~ is the vector bundle \e ~over 
${\rm I} \!{\rm P} ^{n-1} $, the previous vanishing implies in the 
cohomology sequence associated to ~$(*)$~ with $t=0$~ that the hypothesis of the proposition are 
true on ${\rm I}\!{\rm P}^{n-1}$. Hence by induction on $n, 
~~ H^2({\rm I}\!{\rm P}^{n-1} , \, \bigwedge^q E_d(t)_{|{\rm I}\!{\rm P}^{n-1}})=0 
~~ {\rm  for }\break  1\le q\le k ~~ {\rm and} ~~
\forall t\ge 0\,$. From the cohomology sequence associated to $(*) ~ {\rm with} ~ q=k$ ~ we get that the map $H^2(\, {\rm I}\!{\rm P}^n , \, \bigwedge^k E_d(t-1))
\, \lo \, H^2(\, {\rm I}\!{\rm P}^n , \, \bigwedge^k E_d(t)) $~ is surjective
$\forall t \geq 0 $~ and the thesis  follows easily.
\hfill $\diamond$

\vskip1truecm
\noindent
{\it Proof of  Theorem 1.6} \acapo
The implication '$\Longrightarrow$`  is a consequence of 
Proposition 1.5.  To prove the converse, we may apply Proposition 1.7 and then
Proposition 1.5 again. \acapo

\prop {1.8} 
{\it {If ~\fd  ~satisfies ~$N_p$ ~then  ~\fm  ~satisfies ~$N_p ~~\forall m\le n\,$.}}
\pdim 
It follows by the remark of section 2 of [G2] (which is an insight
into representation theory). \hfill $\diamond$

\vfill
\eject

\paragrafo{ 2 - Necessary conditions
 on property {\mathpar ${N_p}$} ~ for the
Veronese embedding {\mathpar $
\varphi _{{\cal O} _{{\rm I}\!{\rm P}^{ n}}(d)}$}}

In this section we will prove the following:

\teo {2.1} {\it {The ~Veronese ~embedding ~\fd ~~ does ~not ~satisfy ~$N_{3d-2} $~ for ~$n\ge 2 \, , ~d\ge 3 ~ \,$.}} 

\pdim
By Proposition 1.8, we can let n=2. By Theorem 1.6 and Serre duality  it is enough to show that  $H^0 \left( \, {\rm I}\!{\rm P} ^2, \, \bigwedge^{K}E_d(d-3) \right) \neq 0 \,$~ with $K:={{d(d-3)}\over 2} \,$. 
 So the theorem will follow from the following:

\lem {2.2} {\it { The bundle ~$\bigwedge ^q E_d (t) $~ has a non  zero
 global section for ~$1\leq q \leq N \, , ~ q+1 \leq {n+t \choose n} ~ and ~ 
t\geq 1\,.$}}

\pdim
The exact sequence  $
0 \, \lo \, \bigwedge ^{q} E_d \, \lo \,
\bigwedge ^{q} S^d V \otimes {\cal O} _{{\rm I}\!{\rm P} ^n}
\, \lo \, \bigwedge ^{q-1} E_d \, (d)
\, \lo \, 0 $~ implies that 
$$
H^0 \left( \bigwedge ^q E_d (t) \right) ~= ~ 
 Ker  \, \left(  \, 
\bigwedge ^q S^d V \, \otimes \, S^t V 
\,  ~~
{
\setbox2=\hbox{$-\hskip -5pt -\hskip -5pt 
\longrightarrow$}
\setbox1=\hbox to \wd2{\hfill $\scriptstyle \alpha _{t}$ 
\hfill} \hbox {\hskip -8pt \vbox {\box1 \vskip -5pt \box2}}
\hskip 2pt} 
 \, 
\bigwedge ^{q-1} S^d V \, \otimes \, S^{t+d} V
 \right)  \, .
$$
Now, there is a Koszul complex
$$
\longrightarrow 
\bigwedge ^{q+1} S^d V \, \otimes \, {\cal O}(t-d) \, \longrightarrow \,
\bigwedge ^q S^d V \, \otimes \, {\cal O}(t)  ~~
{
\setbox2=\hbox{$-\hskip -5pt -\hskip -5pt 
\longrightarrow$}
\setbox1=\hbox to \wd2{\hfill $\scriptstyle a _{t}$ 
\hfill} \hbox {\hskip -8pt \vbox {\box1 \vskip -5pt \box2}}
\hskip 2pt} 
 \, 
\bigwedge ^{q-1} S^d V  \, \otimes \, {\cal O} (t+d)  \,
\longrightarrow
$$
with $\alpha _t \, = \, H^0 (a _t ) \,$.  For $t\geq d \,$, global sections of 
$\bigwedge ^{q+1} S^d V \, \otimes \, {\cal O}(t-d) $~ will therefore give sections of 
$\bigwedge ^q E_d (t) \,$. In particular, for $d=t \,$, we get that for each family 
$s_0 , \cdots \, , s_q $~ of degree $d$~ polynomials,
$$
\sum ^q  _{ i=0}\,(-1)^i ~ s_0 \,\wedge \, \cdots \, \wedge \,  {\widehat {s}_i}\,
\wedge \, \cdots \,\wedge \,s_q  \, \otimes \, s_i
$$
is in the kernel of $\alpha _d \,$.~ Let now be $1 \leq t < d \,$. If we can factor $s_i \, = \, u\, w_i \,$~ with $u$~ of degree $d-t $, ~ then 
$$
\sum ^q  _{ i=0}\,(-1)^i ~ s_0 \,\wedge \, \cdots \, \wedge \,  {\widehat {s}_i}\,
\wedge \, \cdots \,\wedge \,s_q  \, \otimes \, w_i
$$
must be in the kernel of $\alpha _t \,$, and therefore defines a global section of 
$\bigwedge ^q E_d (t) \,$.\break 
Thus to get a non-zero section of 
$\bigwedge ^q E_d (t) \,$, it suffices to find $q+1$~ linearly independent polynomials 
of degree $t \,$, which is possible as soon as 
$q+1 \, \leq \, {{n+t} \choose n}\,.  \hfill \diamond $

\bigskip

\oss {2.3} The bundles $\bigwedge ^{q} E_{d} $~ are semistable
 (see [P], Proposition 5.6), so $H^{0} \left( \bigwedge ^{q} E_{d}
 \, (t) \right) \, = \, 0 ~~ {\rm if}~~ \mu \left( 
\bigwedge ^{q} E_{d} \, (t) \right) \, = \, 
t - {qd \over N} < 0\,.$~ In particular,\acapo
\vskip0.1truecm
\centerline { 
$H^{0} \left( \bigwedge ^{q} E_{d}
 \, (t) \right) \, = \, 0 ~~ \forall t \leq 0 \,.$}
 
\vskip2truecm

\paragrafo{ 3 - Conclusions}

In this section we will fit our results into the literature. \acapo
In particular, we  will prove that:

\teo {3.1} ~ {\it { Let ~$d$~be an integer s.t. ~$d \geq 3 \, . $~ Then the 
Veronese embedding }}\break
{\it { \f2  $\, : \, {\rm I}\!{\rm P} ^2   \,  \hookrightarrow \, $ \pN ~satisfies  property $N_p ~~ if ~and~only~if~~ 0\leq p \leq 3d-3 \,.$\acapo
Moreover, if ~$d=2 \, ,$~the embedding ~ \fpiano2 $\, :
\,  {\rm I}\!{\rm P} ^2 \,  \hookrightarrow  \, $ \p5 ~~satisfies $N_p ~~~\forall p \,.
$}}

\bigskip 
We have the following:

\prop {3.2}  ~~~ (M. Green, C. Birkenhake) \acapo
{\it Let } ~$d \geq 2 ~{\it and } ~~p = \cases  { 3d-3
\cr  2 \cr } ~{\rm if }~ \eqalign { d\geq 3 \cr d=2 \cr} \,.$~
{\it {Then the  complete Veronese embedding ~\f2 $\, : \, 
{\rm I}\!{\rm P} ^2  \; \hookrightarrow  \, $\pN ~~ satisfies property ~$N_p \,$.}} 
\pdim  
See [B1], Corollary $3.2\, $. \acapo
The result follows also applying  
Theorem $3.b.7$~ of [G1](which says that the minimal resolution of
a Veronese variety restricts to the minimal resolution of its curve
hyperplane section) and Theorem $4.a.1$~ of
[G1] (which says that a line bundle of degree $2g+1+p$ on a
curve of genus $g$ satisfies
$N_p$).\hfill $\diamond$

\bigskip
In the same way we get:

\lem{3.3} {\it The Veronese embedding } ~~$\varphi _{{ \cal O}_{{\rm I}\!{\rm P}^3} (3) } \, : \, {\rm I}\!{\rm P}^3 \, \hookrightarrow 
\, {\rm I}\!{\rm P}^{19} $ ~~ {\it satisfies} ~ $N_6$.
\pdim 
The curve hyperplane section of the image of the cubic Veronese 
embedding of ${ {\rm I}\!{\rm P}^{\rm 3}}$ ~is the space curve complete
intersection of two cubics embedded by 
$| {\cal O}_{{\rm I}\!{\rm P}^3} (3)| $  ~ and it  has genus $10$.
The result follows again applying  Theorem $3.b.7 $~ and  Theorem $4.a.1$ ~ of  [G1].\hfill $\diamond$

\lem {3.4} 
{\it {The ideal ~${\cal I}$~ of ~~\fpiano2 $\, ({\rm I}\!{\rm
P} ^2 )$ ~in \p5 ~has the following resolution:}}
$$ O \, \lo \, 
{\cal O}_{ {\rm I}\!{\rm P}^5} (-4) ^{\oplus 3} \, 
\lo \, 
{\cal O}_{ {\rm I}\!{\rm P}^5} (-3) ^{\oplus 8} \, 
\lo \, 
{\cal O}_{ {\rm I}\!{\rm P}^5} (-2) ^{\oplus 6} \, 
\lo \, {\cal I} \, \lo \, 0 \,.
$$
{\it {In particular the Veronese embedding  ~\fpiano2 $\, : \, 
 {\rm I}\!{\rm P}^2 \; \hookrightarrow
 \,$ \p5 ~~satisfies $N_p ~ \forall p \,.$}}
\pdim 
Easy computation. \hfill $\diamond$ 
\bigskip

\bigskip
\noindent
{\it Proof ~of ~Theorem ~3.1} \acapo
By Proposition $3.2$~ and Lemma $3.4$~ we just need to show
that if $d\geq 3$~ then property $N_{p}$~ does not hold 
for $p \geq 3d-2 \,.$~  But this is exactly the bound coming
from Theorem $2.1\, . \hfill \diamond$

\vskip1truecm

When $d=2$~ the minimal free resolution of the quadratic Veronese variety is known from the work of Jozefiak, Pragacz and Weyman [JPW],
in which they prove a conjecture made by Lascoux. As a  corollary of the above paper we have  the following result
(which agrees with our Conjecture formulated in the Introduction):

\teo {3.5} {\it {The quadratic Veronese embedding }} ~~
$\varphi _{{\cal O} _{{\rm I}\!{\rm P}^n  (2) }} \, : \, $ \pn $\, $ \hr
$\, $\pN ~~ {\it satisfies} 
~~$N_p ~~  {\it if ~ and ~ only ~ if } ~~p \leq 5 ~  {\it when}  ~n\geq 3
~~~{\it and} ~~ \forall p ~{\it when}~ n=2 \,.
\hfill \diamond $

\bigskip
The following nice characterization, probably well known, 
was found during discussions with E. Arrondo: 

\teo {3.6} 
{\it {The only (smooth)  varieties in \pn ~ such that ~$N_p$ ~holds for every ~$p\ge 0$ ~ are the quadrics, the rational normal scrolls and the  Veronese surface in ~ \p5 .}}
\pdim
Suppose $X$~ is a variety satisfying $N_p$ ~ for every $p \ge 0 \,$.  
Then 
$H^i \left( {\cal O}_X (t) \right) =0 ~{\rm for} \break
t\ge 0 ~{\rm and}~ 
1\le i \le dim\, X -1 \,$. ~ Hence from Theorem $3.b.7$~ in [G1]  it follows
that the minimal free resolution of $X$ ~ restricts to the minimal resolution
of its generic curve section $C \,$.  This implies that $H^1 \left( {\cal O}_C
\, \right) =0$~ and $C$~ is linearly normal, hence $C$~ is a rational normal
curve. In particular, $X$~ has minimal degree and we get the result. $\hfill
\diamond $

\medskip
We remark that the only Veronese varieties appearing in Theorem $3.6$~ are the rational normal curves and the Veronese surface in \p5 $\,$.

\vskip2truecm
\paragrafo {References} {}

\parskip 15pt plus 8pt minus 3pt
 
\cite{B1} 
C. Birkenhake  \hfill\break  
{\it  Linear systems on projective spaces}\hfill\break  
Manuscripta Math. 88 (1995),  177-184 

\cite{B2}
C. Birkenhake  \hfill \break 
{\it Noncomplete linear systems on abelian varieties}
\hfill \break
Trans. A.M.S.  vol 348, Number 5  (1996),  1885-1908 

\cite{EL}
L. Ein, R. Lazarsfeld  \hfill \break 
{\it Syzygies and Koszul cohomology of smooth projective varieties of arbitrary 
dimension}
\hfill \break
Invent. Math. 111 (1993), 51-67

\cite{G1}
M. Green  \hfill \break 
{\it Koszul cohomology and the geometry of projective varieties }
\hfill \break 
J. Diff. Geometry 19 (1984),  125-171.

\cite{G2}
M. Green  \hfill \break 
{\it Koszul cohomology and the geometry of projective varieties, II }
\hfill \break 
J. Diff. Geometry 20 (1984),  279-289.

\cite{GL}
M.Green, R.Lazarsfeld \hfill \break
{\it Some results on the syzygies of finite sets and algebraic curves} 
\hfill\break
Compositio Math. 67 (1988), 301-314.

\cite{JPW}
T. Jozefiak, P. Pragacz, J. Weyman \hfill \break
{\it Resolutions of determinantal varieties and tensor complexes
associated with symmetric and antisymmetric matrices }\hfill \break
Asterisque 87-88 (1981), 109-189.

\cite{M}
L. Manivel \hfill \break
{\it On the syzygies of flag manifolds}\hfill \break
Proc. AMS 124 (1996), 2293-2299

\cite{P} R. Paoletti \hfill \break  
{\it Stability of a class of homogeneous vector bundles 
on ${\rm I}\!{\rm P}^n $} \hfill \break 
Bollettino U.M.I. (7) 9-A (1995), 329-343

\vskip2truecm 

{\bf Authors' addresses}
\vskip0.5truecm 
\noindent
\hbox  {\hsize15truecm
\vbox {
\hsize6truecm
Giorgio Ottaviani\hfill \break  
Dipartimento di Matematica \hfill \break 
"U. Dini" \hfill \break
Universit\`a di Firenze\hfill \break 
viale Morgagni 67/A\hfill \break 
I 50134 FIRENZE\hfill \break  
ottavian@math.unifi.it
} 
\hfill 
\vbox  {
\hsize6truecm
Raffaella Paoletti\hfill \break 
Dipartimento di Matematica ed Appl. per l'Architettura \hfill \break 
Universit\`a di Firenze\hfill \break 
piazza Ghiberti\hfill \break 
50122 Firenze \hfill \break 
raffy@math.unifi.it}
}
\vfill
\eject
\bye